\documentclass{article}
\newcommand{\proof}{\noindent {\bf Proof: }}
\newcommand{\note}{\noindent {\bf Notation: }}
\newcommand{\remark}{\noindent {\bf Remark:}}
\newcommand{\corollary}{\noindent {\bf Corollary:}}
\newcommand{\ws}{\hspace{4pt}}
\newtheorem{theorem}{Theorem}

\newtheorem{lemma}{Lemma}
\newtheorem{defi}{Definition}
\usepackage[]{amssymb}

\begin{document}

\title{M\"untz-type theorems on the half-line with weights}
\author{\'Agota P. Horv\'ath \footnote{Supported by Hungarian National Foundation for Scientific Research, Grant No. K-61908.\newline Key words:M\"untz theorem, complete system, weighted spaces. \newline  2000 MS Classification: 42A65, 42A55}}\date{}
\maketitle
\begin{center}\small{Department of Analysis, Budapest University of Technology and Economics}\end{center}
\begin{center}\small{H-1521 Budapest, Hungary  }\end{center}
 \begin{center}\small{e-mail: ahorvath@renyi.hu}\end{center}

\begin{abstract} We consider the linear span $S$ of the functions $t^{a_k}$ (with some $a_k>0$) in weighted $L^2$ spaces, with rather general weights. We give one necessary and one sufficient condition for $S$ to be dense. Some comparisons are also made between the new results and those that can be deduced from older ones in the literature.
\end{abstract}

\section{Introduction}

The first "if and only if" solution to a problem of S. N. Bernstein \cite{be} was given by Ch. H. M\"untz \cite{mu}:

\medskip

\noindent{\bf Theorem A}

\noindent{\it Let $0=\lambda_0<\lambda_1<, \dots$ be an increasing sequence of real numbers. The linear subspace $\mbox{span} \{t^{\lambda_k} : k=0,1,\dots\}$ is dense in $C([0,1])$, if and only if $\sum_{k=1}^{\infty}\frac{1}{\lambda_k} = \infty$.}

\medskip

This classical result was first proved in $L_2[0,1]$ and then extended to $C[0,1]$, as stated above. Also, it was stated only for increasing sequences $\lambda_k$. Subsequently, this theorem has had several different proofs and generalizations, and there are several surveys in this topic (see for instance the papers of J. Almira and A. Pinkus \cite{a}, \cite{p}).

\medskip

On $C[0,1]$ and $L_p(0,1)$, "full M\"untz theorems", i.e. theorems with rather general exponents, were later proved by eg. P. Borwein, T. Erd\'elyi, W. B. Johnson and V. Operstein (\cite{be1}, \cite{ej}, \cite{e}, \cite{o}). Versions of M\"untz's theorem on compact subsets of positive measure \cite{be2}, \cite{be3}, and on countable compact sets \cite{a1} were also proved. Further results can be found for instance in the monographs of P. Borwein, T. Erd\'elyi \cite{be4}, and B. N. Khabibullin \cite{kh}.

\medskip

In this paper we are interested in M\"untz-type theorems on $(0,\infty)$. Several papers were written in the '40s on the completeness of the set $\{t^{\lambda_k}e^{-t}\}$ in $L_2(0,\infty)$ (see eg. \cite{fu}, \cite{bo}, \cite{bopo}). In particular, we will use some ideas of W. Fuchs. His theorem is the following:

\medskip

\noindent{\bf Theorem B}

{\it Let $a_k$ be positive numbers, such that $a_{k+1}-a_k >c>0$ $(k=1,2,\dots),$ and let $\log \Psi(r) = 2\sum_{a_k<r}\frac{1}{a_k},$ if $r>a_1$, and $\log \Psi(r) = \frac{2}{a_1}$ if $r\leq a_1$. Then $\{e^{-t}t^{a_k}\}$ is complete in $L_2(0,\infty)$, if and only if
$$\int_1^{\infty}\frac{\Psi(r)}{r^2} = \infty$$}

A. F. Leontev \cite{le} and G. V. Badalyan \cite{ba} proved similar theorems with more general weights (the weight being $e^{-t}$ in the above theorem). In 1980, by the Hahn-Banach theorem technique, R. A. Zalik \cite{za} proved a M\"untz type theorem on the half-line with weights $|w| \leq c\exp(-|\log t|^a)$ $(a>0)$. In 1996 Kro\'o and Szabados \cite{ks} also had a related result on $(0, \infty)$.

\medskip

In Theorem \ref{thm1} and Theorem \ref{thm2} below we will prove M\"untz-type theorems on the half-line with more general weights, which generalize all the results mentioned above.

\medskip

Closely related to our topic (by a $\log t$ substitution) are the results on the whole real line for exponential systems. The basic paper in this respect was written by P. Malliavin \cite{ma}, and by this tool there are some nice generalizations of the above mentioned results, for instance by B. V. Vinnitskii, A. V. Shapovalovskii \cite{visha}, by G. T. Deng \cite{de}, and by E. Zikkos \cite{zi}.

\section{Definitions, Results}

Let us begin with a rather general definition. Some specific examples are given subsequently.

\begin{defi}We say that a weight function $w(t)=\nu(t)\mu(t)$ is admissible on $[0,\infty)$, if $\nu(t)$ and $\mu(t)$ are positive and continuous on $(0,\infty)$, $w^2$ has finite moments,
and there is a function $\gamma$ on $[0,\infty)$, such that
$$\gamma (t) = \sum_{k=0}^{\infty}c_kt^{\gamma_k},$$
where $c_k >0$ for all $k$, and $ 0\leq \gamma_0<\gamma_1 <\gamma_2 < \dots$, and there is a $C_0>0$ such that $\forall t>C_0$
\begin{equation}\frac{1}{w^2(t)} \leq \gamma(t)\end{equation}
and there is a $C>1$, such that
\begin{equation}\int_0^{\infty}\gamma\left(\frac{t}{C}\right)w^2(t)dt <\infty. \end{equation}
Furthermore we require that
\begin{equation} \lim_{t\to 0+}\mu(t) \in (0,\infty), \end{equation}
and there is an $a>0$ such that
\begin{equation}\int_0^1\left(\frac{t^{a-1}}{\nu(t)}\right)^2 < \infty. \end{equation}
\end{defi}

\medskip

Here and in the followings $C, C_i$ and $c$ are absolute constants, and the value of them will not be the same at each occurrence.

\medskip

\remark

If $\nu(t)\equiv 1$ (as in Theorem B) then we can choose $a=1$. Also, it is easy to see that we can always assume $a\geq 1.$

\medskip

\noindent {\bf Examples :}

$w(t)= t^{\beta}e^{-Dt^{\alpha}}$, where $\beta > -\frac{1}{2}$ and $\alpha > 0$ is admissible, namely it has finite moments, and $\gamma(t)=e^{3Dt^{\alpha}}$ serves the purpose. When $\beta =0$ and $D=\alpha =1$, we get back the original case of Fuchs (Theorem B). When $\beta > -\frac{1}{2}$, $D=\frac{1}{2}$ and $\alpha =1$, then $w^2$ is a Laguerre weight. When $\beta =0$, $D>0$ and $\alpha > 1$, then $w$ is a Freud weight.

Let $w(t) =(4+\sin t)t^{\beta}\prod_{k=1}^ne^{-D_kt^{\alpha_k}}$, and let us assume, that $\beta > -\frac{1}{2}$,  and $0\leq \alpha_1< \alpha_2 \dots < \alpha_n$, and $D_n>0$.
Then $w$ is admissible, and $e^{Dt^{\alpha_n}}$ is a suitable choice for $\gamma(t)$, if $D$ is large enough. In particular, if $w(t)=t(4+\sin t)e^{-t}$ then the second derivative of $-\log(w(e^{t}))$ takes some negative values on $(A, \infty)$ for any $A>0$. This means that the results of \cite{zi} are not applicable in this case.

\medskip

\begin{defi}
Let $w$ be a positive continuous weight function with $w^2$ having finite moments. Then define $\varphi(x)$ and $K(x)$ corresponding to $w(x)$ as

\begin{equation}\varphi(x) = \left(\int_0^{\infty}t^{2x}w^2(t)dt\right)^{\frac{1}{2x}} = \left(K(x)\right)^{\frac{1}{2x}}, \ws \ws x > 0.\end{equation}\end{defi}

\medskip

Furthermore let us define another property of a weight function. The classical weight functions, and also our examples above, fulfil this "normality" condition, as we can see later.

\medskip

\begin{defi} Let us call a weight function $w^2$ with finite moments "normal", if the largest zero of the $n^{th}$ orthogonal polynomial $(x_{1,n})$ with respect to $w^2$, can be estimated as:
$$x_{1,n} \leq e^{cn},$$
where $c=c(w)$ is a positive constant independent of $n$. \end{defi}

\medskip

\remark
\medskip

In the cases of Laguerre and Freud weights $x_{1,n} \leq cn^{\lambda}$, where $\lambda=\lambda(w)$ is a positive constant depending on the weight function, moreover the same estimation is valid for a more general classes of weights on the real axis (\cite{lelu} p. 313. Th. 11.1). As an application of the result of A. Markov (\cite{sze} p. 115. Th. 6.12.2), we can get a similar estimation for the examples above; for instance $w(x)=x^{\gamma}e^{-x^{\alpha}}$, there is a $\beta >0$ such that with $W(x)= x^{\beta}e^{-x}$, the quotient $\frac{W}{w}$ is increasing on $(0, \infty)$; if $w(x)= x(4+\sin x)e^{-x}$, then the corresponding $W$ can be $W(x) = x^2e^{-\frac{x}{2}}$.

\medskip

\begin{defi}Let $\{a_k\}_{k=1}^{\infty}$ be positive numbers in increasing order. We define (as in \cite{fu} and Theorem B above)
\begin{equation}m(r) = \left\{\begin{array}{ll}\frac{1}{a_1}, \ws \mbox{if} \ws 0 \leq r \leq a_1\\
\sum_{a_k<r}\frac{1}{a_k}, \ws \mbox{if} \ws r>a_1 \end{array}\right.\end{equation}
and let
\begin{equation} \Psi(r) = e^{2m(r)}.\end{equation}

\end{defi}

Let us also introduce the following notations:

\note

Let $w$ be a positive continuous weight function, and let us define the weighted $L^2$ space as $L^2_w(0,\infty) = \{f | fw \in L^2(0,\infty)\}$, and $\|f\|_{2,w}=\|fw\|_{2, (0,\infty)}.$

\medskip

\begin{equation}S=\mbox{span}\{t^{a_k}: k=1,2,\dots\}\end{equation}
with $0<a_1 < a_2 < \dots$.

\medskip

We are now in position to state the main results of this note (the proofs will follow in the next Section).

\begin{theorem}\label{thm1}Let $w$ be an admissible and normal weight function on $[0,\infty)$.
If there exists a monotone increasing function $f$ on $[0,\infty)$, such that for all $0< x \leq r$
\begin{equation}x\log\frac{\Psi(r)}{\varphi(x)} \leq f(r),\end{equation}
and
\begin{equation}\int_1^{\infty}\frac{f(r)}{r^2} < \infty , \end{equation}
then $S$ is incomplete in $L^2_w(0,\infty)$. \end{theorem}

\medskip

This result is then nicely complemented by the following positive result.

\medskip

\begin{theorem}\label{thm2}Let $w$ be positive and continuous on $(0,\infty)$, such that $w^2$ has finite moments. Let us suppose that there is a constant $d>0$ such that
\begin{equation} a_{k+1} - a_k >d \end{equation}
If there exists a monotone increasing function $h$ on $[0,\infty)$, for which
\begin{equation} C < \frac{h(r)}{h(r_1)} < D,
\ws\ws\ws \mbox{for}\ws\ws\ws  \frac{1}{2} \leq \frac{r}{r_1} \leq 2 \end{equation}
with some $0< C,D$, and there are $\alpha, C, c >0$, such that for all $0< x \leq r$
\begin{equation}0< h(r) \leq C^{\frac{1}{x}}\frac{cx}{\varphi^{\alpha}(x)} \Psi^{\alpha}(r),\end{equation}
and
\begin{equation}\int_1^{\infty}\frac{h(r)}{r^2} = \infty , \end{equation}
then $S$ is complete in $L^2_w(0,\infty)$. \end{theorem}

\medskip

Comparing the conditions of the above theorems we conclude the following:

\medskip

\corollary

\noindent If $w$ is admissible and normal on $(0,\infty)$, and there is a $d$ such that $a_{k+1} - a_k >d >0$, and $f(r) = ch(r)$, where $h$ has the same properties as in Theorem 2, then $S$ is dense in $L^2_w(0,\infty)$ if and only if $\int_1^{\infty}\frac{h(r)}{r^2} = \infty $.

\medskip

\medskip

\remark

\medskip
\noindent (1) Let
$$B_{\alpha}(r) = \inf_{x\in (0,r)}C^{\frac{1}{x}}\frac{x}{\varphi^{\alpha}(x)}. $$
Then assuming (11) and (12), if there exists a $0\leq h(r) \leq c B_{\alpha}(r)\Psi^{\alpha}(r)$, for which (14) is valid, then $S$ is dense in $L^2_w(0,\infty)$.

\medskip

\noindent (2) Theorem 2 can be stated also in $L^p_w(0,\infty)$, with $1\leq p < \infty$, and in $C_{w, (0,\infty)}$ with the same proof. That is,
let us define
$$\varphi_p(x) = \left(\int_0^{\infty}t^{px}w^p(t)dt\right)^{\frac{1}{px}}, \ws \ws x > 0, \ws\ws 1\leq p < \infty$$
and
$$\varphi_c(x) = \left(\sup_{t>0}t^{x}w(t)dt\right)^{\frac{1}{x}}, \ws \ws x > 0.$$
Using the standard the notations $L^p_w(0,\infty)=\{f: \|fw\|_{p,(0,\infty)}<\infty\}$, and  $C_{w,(0,\infty)}=\{f \in C(0,\infty): \lim_{t \to 0+ \atop t \to \infty} f(t)w(t) =0\}$, we can formulate the following theorem:

\begin{theorem}Let $w$ be positive and continuous on $(0,\infty)$, and let us assume that $t^xw(t) \in L^p_{(0,\infty)}$ in the $L^p_w$-case, and that for all $a>0$ $\lim_{t \to 0+ \atop t \to \infty} t^aw(t) =0$ in the $C_w$-case. Furthermore let $\{a_k\}$ be a sequence of positive numbers for which (11) is satisfied. If there is a monotone increasing function $h$ on $(0,\infty)$ with the properties (12) and (14), and for which there are $\alpha, C, c >0$, such that for all $0< x \leq r$
$$0< h(r) \leq C^{\frac{1}{x}}\frac{cx}{\varphi_{p/c}^{\alpha}(x)} \Psi^{\alpha}(r),$$
then $S$ is complete in $L^p_w(0,\infty)$/in $C_{w,(0,\infty)}$.
\end{theorem}

\medskip

\noindent (3) If $B_{\alpha}(r)>B>0$ ($\forall r\geq 1$), then $h(r) = B\Psi^{\alpha}(r)$. This is the situation when $w(t) = e^{-Dt^{\alpha}}.$ Furthermore with suitable $D$ and $\alpha$ $\inf_{x\in (0,r)}\frac{x}{\varphi^{\alpha}(x)}>B>0$.
In this case
$$K(x) = \int_0^{\infty}t^{2x}e^{-2Dt^{\alpha}}dt = \frac{1}{\alpha(2D)^{\frac{2x+1}{\alpha}}}\Gamma\left(\frac{2x+1}{\alpha}\right)$$
By Stirling's formula (see eg. \cite{he})
$$\frac{x}{\varphi^{\alpha}(x)}= \frac{x}{\left( \frac{\sqrt{2\pi}\left(\frac{2x+1}{\alpha}\right)^{\frac{2x+1}{\alpha}-\frac{1}{2}}e^{- \frac{2x+1}{\alpha}}e^{J\left(\frac{2x+1}{\alpha}\right)}}{\alpha(2D)^{\frac{2x+1}{\alpha}}}\right)^{\frac{\alpha}{2x}}}= (*),$$
where $J$ is the Binet function. For $x>0$ we have $ 0<J(x) < \frac{1}{12x}$. That is,
$$(*) \geq  \frac{2De\alpha x}{2x+1}\left(\frac{2De\alpha}{2x+1}\left(\frac{\alpha(2x+1)}{2\pi} \right)^{\frac{\alpha}{2}}\frac{1}{e^{\frac{\alpha^2}{12(2x+1)}}}\right)^{\frac{1}{2x}}= b(D, \alpha, x)$$
and $b(D, \alpha, x)$ tends to $De\alpha$ when $x$ tends to infinity, and if \\ $C(\alpha, D)=\frac{\sqrt{2De\alpha}}{e^{\frac{\alpha^2}{12}}\left(\frac{2\pi}{\alpha}\right)^{\frac{\alpha}{4}}}>1$ then $\lim_{x\to 0+}\frac{x}{\varphi^{\alpha}(x)} =\infty$. (In the case of Fuchs, $\alpha=D=1$, $C(\alpha, D)>1$.)

\medskip

\noindent (4) For $\alpha=D=1$, $h(r)=f(r)= \Psi(r), r\geq 0.$ By the substitution $t=Du^{\alpha}$ (without any further restrictions on the exponents $a_k$ for $\alpha \geq 1$, and with the restriction $a_k\neq \frac{1}{2}(\frac{1}{\alpha}-1)$ for $0<\alpha <1$), after some obvious estimations one can deduce from the result of Fuchs (Theorem B), that $\{t^{a_k}e^{-Dt^{\alpha}}\}$ is dense if and only if $\int_1^{\infty}\frac{\Psi^{\alpha}(r)}{r^2} =\infty.$ We get the same from Theorems 1 and 2. After the third remark we need to check the assumptions of Theorem 1. Now $\varphi^x(x)= \sqrt{K(x)}\leq (cx)^{\frac{x}{\alpha}}$, and so
$$\left(\frac{\varphi(x)}{\Psi(r)}\right)^x \leq \left(\frac{cx}{\Psi^{\alpha}(r)}\right)^{\frac{x}{\alpha}}.$$

\medskip

\begin{theorem}\label{thm4}
With the notations of Theorem 1
\begin{equation} f(r)=C+r\max\left\{\frac{1}{2}\frac{K^{'}}{K}(r), 2m(r)\right\}-\frac{1}{2}\log K(r)\end{equation}
is a good choice for $f(r)$ with a suitable $C$. That is, if $w$ is admissible on $[0,\infty)$, and
\begin{equation}\int_1^{\infty}\frac{r\max\left\{\frac{1}{2}\frac{K^{'}}{K}(r), 2m(r)\right\}-\frac{1}{2}\log K(r)}{r^2} <\infty\end{equation}
then $S$ is incomplete in $L^2_w(0,\infty)$.
\end{theorem}

\medskip

\remark

\medskip

If $w(t) = e^{-Dt^{\alpha}}$, and $\frac{1}{2}\frac{K^{'}}{K}(r) > 2m(r)$ on a set $H$, then on $H$
$$f(r) = r\frac{1}{2}\frac{K^{'}}{K}(r) -\frac{1}{2}\log K(r)$$ $$= \frac{r}{\alpha}\left(\frac{\Gamma^{'}}{\Gamma}\left(\frac{2r+1}{\alpha}\right)-\log(2D)\right)-\frac{1}{2}\log \left(\frac{1}{\alpha(2D)^{\frac{2x+1}{\alpha}}}\Gamma\left(\frac{2x+1}{\alpha}\right)\right)$$
$$
=\frac{r}{\alpha}\log\frac{2r+1}{\alpha}-\frac{r}{2(2r+1)}-\frac{r}{\alpha}I\left(\frac{2r+1}{\alpha}\right) -\frac{r}{\alpha}\log(2D)
$$ $$-\frac{1}{2}\log\frac{\sqrt{2\pi}\left(\frac{2x+1}{\alpha}\right)^{\frac{2x+1}{\alpha}-\frac{1}{2}}e^{- \frac{2x+1}{\alpha}}e^{J\left(\frac{2x+1}{\alpha}\right)}}{\alpha(2D)^{\frac{2x+1}{\alpha}}}$$
$$ =\frac{r}{\alpha} - \frac{2-\alpha}{4\alpha}\log\frac{2r+1}{\alpha} + O(1)$$
(In the last step we used that $\frac{\Gamma^{'}}{\Gamma}(z) = \log z-\frac{1}{2z} -I(z)$ where $I(z)= \\ 2\int_0^{\infty}\frac{t}{(t^2+z^2)(e^{2\pi t}-1)}dt$ (see eg \cite{wiwa}).) That is, if $H$ is large then the integral in (10) is divergent.

\section{Proofs}

For the proof of the first theorem, at first we need a lemma:

\medskip

\begin{lemma} Let $a=m$ be a positive integer. If $w^2$ is a continuous, positive, normal weight function on $(0,\infty)$ with finite moments, then there is a function $b(z)$ such that $\frac{1}{b(z)}$ is regular on $\Re z \geq -a$, and it fulfils the inequality on $\Re z \geq -\frac{1}{2}$:
$$\sqrt{\frac{K(x+a)}{K(x)}} \leq |b(z)|,$$
where $z=x+iy$.
\end{lemma}

\medskip

\proof

At first let $x=n$ be also a positive integer. Then, using the Gaussian quadrature formula on the zeros of the $N^{th}$ orthogonal polynomials ($x_{1,N} > \dots >x_{k,N}> \dots > x_{N,N}$) with respect to $w^2$, where $N=n+m+1$, we get, that
$$\frac{K(n+m)}{K(n)}= \frac{\int_0^{\infty}t^{2(n+m)}w^2}{\int_0^{\infty}t^{2n}w^2}= \frac{\sum_{k=1}^N\lambda_{k,N}x_{k,N}^{2(n+m)}}{\sum_{k=1}^N\lambda_{k,N}x_{k,N}^{2n}}\leq x_{1,N}^{2m},$$
that is, by the condition of "normality"
$$\sqrt{\frac{K(n+m)}{K(n)}} \leq e^{cNm}.$$
Now we can consider, that
$\frac{K(x+a)}{K(x)}$ is increasing on $\Re z > -\frac{1}{2}$, namely
$$\left(\frac{K(x+a)}{K(x)}\right)^{'} = \frac{K(x+a)}{K(x)}\left(\frac{K^{'}}{K}(x+a)-\frac{K^{'}}{K}(x)\right),$$
which is nonnegative, because $\frac{K^{'}}{K}$ is increasing. The last statement can be seen by the Cauchy-Schwarz inequality, that is the derivative of $\frac{K^{'}}{K}$ is nonnegative:
$$\left(2\int_0^{\infty}t^{2x}|\log t|w^2(t)dt\right)^2 \leq \int_0^{\infty}t^{2x}w^2(t)dt\int_0^{\infty}t^{2x}4\log^2 tw^2(t)dt.$$
So with $a=m$ and $x>0$, 
$$\sqrt{\frac{K(x+a)}{K(x)}} \leq e^{ca(a+1+\lceil x\rceil)} \leq e^{ca(a+2+x)}= C(a)\left|e^{caz}\right|.$$

\medskip

\remark

\medskip

\noindent (1) If $\left(\frac{\varphi(2x)}{\varphi(x)}\right)^x$ does not grow too quickly, then one can choose $b(z)= c(a)b_1(z),$ where $b_1(z)$ is independent of $a$, because
\newline
$\sqrt{\frac{K(x+a)}{K(x)}} \leq K^{\frac{1}{4}}(2a)\frac{K^{\frac{1}{4}}(2x)}{K^{\frac{1}{2}}(x)}=c(a)\left(\frac{\varphi(2x)}{\varphi(x)}\right)^x$

\medskip

\noindent (2) Usually we can give $b(z)$ in polynomial form, for instance if $w(t) = e^{-Dt^{\alpha}}$ then $\sqrt{\frac{K(x+a)}{K(x)}} = \frac{1}{(2D)^{\frac{a}{\alpha}}}\sqrt{\frac{\Gamma\left(\frac{2x+1}{\alpha}+\frac{2a}{\alpha}\right)}{\Gamma\left(\frac{2x+1}{\alpha}\right)}}\leq c(2x+1+2a)^{\frac{a}{\alpha}}$, and so $b(z)= c(z+2a)^n$, where $n>\frac{a}{\alpha}$ is an integer.

\medskip

\proof  of Theorem 1.

\medskip

Let us extend $f(r)$ to $\mathbb{R}$ as $f(-r) = f(r)$. Let $a\geq 1$ be as in (4). Furthermore let $a$ be an integer. Because $\int_1^{\infty}\frac{f(r)}{r^2} < \infty ,$ the function
\begin{equation} p(z) = p(x+iy)= p(re^{i\vartheta})=\frac{2}{\pi}(x+a)\int_{-\infty}^{\infty}\frac{f(t)}{(x+a)^2 + (t-y)^2}dt\end{equation}
is harmonic on $\Re z > -a$. Since $f(t)$ is increasing, and $x^2+y^2=r^2$
$$p(z) \geq \frac{2}{\pi} f(r) \int_{|t|>r}\frac{x+a}{(x+a)^2 + (t-y)^2}dt $$ $$= f(r)\frac{2}{\pi}\left(\pi - \left(\arctan\frac{r-y}{x+a} + \arctan\frac{r+y}{x+a}\right)\right) > f(r).$$
(In the last inequality we applied the height theorem of a triangle.)
Let us choose $q$ so that $-p+iq$, and hence $g(z)=g_a(z) = e^{-p+iq}$, be regular on $\Re z> -a$. According to the assumptions of Theorem 1, for this $g(z)\not\equiv 0$ on $\Re z \geq -a$ we have that
\begin{equation} |g(z)| \leq e^{-f(r)} \leq \left(\frac{\varphi(x)}{\Psi(r)}\right)^x \ws\ws\ws \Re z\geq 0.\end{equation}
We will show that in this case $S$ is not dense. For this let us define a regular function on the half plane $\Re z\geq 0$ by
\begin{equation}H(z) = \prod_{k=1}^{\infty}\frac{a_k-z}{a_k+z}e^{\frac{2z}{a_k}}\end{equation}
According to a Lemma of Fuchs (\cite{fu} L.5)
\begin{equation}|H(z)|\leq (C\Psi(r))^x \ws\ws\mbox{on} \ws\ws \Re z\geq 0. \end{equation}
Let us replace the $a_k$-s in the definition of $H(z)$ by $a_k+a$, and let us denote the new function by $H^*(z)$. Now, with the help of $g$ and $H^*$ we can define a function $G(z)=G_a(z)$ which is regular on $\Re z \geq -a$:
\begin{equation} G(z)= \frac{g(z+a)H^*(z+a)}{b(z)C_1^{z+a}}, \end{equation}
where, according to Lemma 1, $\frac{1}{b(z)}$ is regular on $\Re z \geq -a$, and on $\Re z \geq -\frac{1}{2}$ we have
\begin{equation}\sqrt{\frac{K(x+a)}{K(x)}} \leq |b(z)|\end{equation}
Because for an $a>0$ $\frac{K(x+a)}{K(x)}$ is positive, and it tends to zero, when $x$ tends to $-\frac{1}{2}$, according to Lemma 1,
we can suppose that $|b(z)| > \delta >0$ on $\Re z\in \left[-a, -\frac{1}{2}\right]$.

Now, because $|H^*(z+a)| \leq (C\Psi(r))^{x+a}$ $(x\geq -a)$ (see (20)), we have that if $C_1$ is large enough, than according to (22)
\begin{equation} |G(z)| \leq (\varphi(x))^x \ws\ws\ws\mbox{on}\ws\ws \Re z > -\frac{1}{2} ,\end{equation}
and because $ a> \frac{1}{2}$, on $\Re z \in \left[-a, - \frac{1}{2}\right]$ :
\begin{equation}|G(z)| \leq \frac{\left(\varphi(x+a)\right)^{x+a}}{|b(z)|}\leq \frac{1}{\delta}\max_{x \in \left[-a, - \frac{1}{2}\right]}\sqrt{K(x+a)}= M\end{equation}

In the followings we will show that if there exists a function $G$ which is not identically zero, and is regular on $\Re z \geq -a$, and fulfils the equations $G(a_k)=0$ $(k=1,2,\dots)$, and the inequalities (23) and (24) are valid, then $S$ is not complete.

For the purpose of showing this, we need to construct a function $0\not\equiv k(t) \in L^2_w(0,\infty)$ such that $\int_0^{\infty}t^{a_k}k(t)w^2(t)=0$ for $k=1,2,\dots$. We give $k(t)$ by the inversion formula for the Mellin transform of $\frac{G(z)}{(1+a+z)^2}$: on $\Re z \geq -a$ let us define the function $u(t)$ by an integral along a line parallel with the imaginary axis
\begin{equation}t\nu(t)u(t)= \frac{1}{2\pi}\int_{-\infty}^{\infty}\frac{G(z)}{(1+a+z)^2}t^{-z}dy\end{equation}
It  can be easily seen (by taking the integral round a rectangle $x_k \pm iL$ $k=1,2$, where $L \to \infty$) that the integral is independent of $x$. Let us choose
\begin{equation}k(t) = \frac{\nu(Ct)u(Ct)}{w^2(t)},\end{equation}
where $C$ is the same as in  (2). Using that
$$\frac{G(z)}{(1+a+z)^2} = \int_0^{\infty}\nu(t)u(t)t^zdt,$$ we have that
$$\int_0^{\infty}t^{a_k}k(t)w^2(t)dt= \frac{1}{C_1^{a_k+a}}\int_0^{\infty}v^{a_k-1}vu(v)\nu(v)dv$$ \begin{equation}=\frac{1}{C_1^{a_k+a}}\frac{G(a_k)}{(1+a+a_k)^2}=0\end{equation}
We have to show, that $k(t) \in L^2_w(0,\infty)$.
\begin{equation}\|k\|^2_{2,w}=\int_0^{\infty}\frac{u^2(Ct)\nu^2(Ct)}{w^2(t)}dt = \int_0^{\frac{A}{C}}(\cdot) + \int_{\frac{A}{C}}^{\infty}(\cdot)= I + II, \end{equation}where $A = \max\{1, CC_0\}$.

According to (1), and by the positivity of the coefficients in $\gamma$,
$$II \leq \int_{\frac{A}{C}}^{\infty}\nu^2(Ct)u^2(Ct)\sum_{k=0}^{\infty}c_k t^{\gamma_k}dt \leq \sum_{k=0}^{\infty}\frac{c_k}{C^{\gamma_k+1}} \int_{A}^{\infty}t^{\gamma_k}\nu^2(t)u^2(t)dt$$
\begin{equation}=\sum_{k \atop \gamma_k<\frac{1}{3}}(\cdot) + \sum_{k \atop \gamma_k\geq \frac{1}{3}}(\cdot)= S_1 + S_2\end{equation}

Using Parseval's formula for the Mellin transform (see e.g. \cite{he})
$$\int_0^{\infty}t^{2x+1}\nu^2(t)u^2(t)=\frac{1}{2\pi}\int_{-\infty}^{\infty}\left|\frac{G(z)}{(1+a+z)^2}\right|^2dy $$ \begin{equation}\leq c\left(\varphi(x)\right)^{2x} \leq c\left(\varphi(x+\frac{1}{2})\right)^{2x+1},\end{equation}
where the equality is valid on $\Re z \geq -a$, the first inequality is on $\Re z > -\frac {1}{2}$, and the last inequality is on $\Re z \geq -\frac {1}{3}$ say, where we used again that $\frac{K(x+a)}{K(x)}$ is increasing, that is
$$0<c \leq \frac{K(\frac{1}{6})}{K(-\frac{1}{3})}\leq \frac{K(x+\frac{1}{2})}{K(x)}$$
Therefore, by (30)

$$S_2\leq c\sum _{\gamma_k \atop k\geq \frac{1}{3}}\frac{c_k}{C^{\gamma_{k+1}}}\int_0^{\infty}t^{\gamma_k}\nu^2(t)u^2(t)dt \leq c \sum_{k=0}^{\infty} \frac{c_k}{C^{\gamma_{k+1}}}
\left(\varphi\left(\frac{\gamma_k}{2}\right)\right)^{\gamma_k}$$ $$\leq c\sum _{k=0}^{\infty}\frac{c_k}{C^{\gamma_{k+1}}}\int_0^{\infty}t^{\gamma_k}w^2(t)\leq c \sum _{k=0}^{\infty}c_k\int_0^{\infty}\left(\frac{t}{C}\right)^{\gamma_k}w^2(t)$$ \begin{equation}=c\int_0^{\infty}\gamma\left(\frac{t}{C}\right)w^2(t)<\infty\end{equation}
To estimate $S_1$ and $I$, we will use that by (25), with $x = -\frac{1}{3}$
\begin{equation}\nu^2(t)u^2(t)\leq c t^{-\frac{4}{3}}\varphi\left(-\frac{1}{3}\right)^{-\frac{2}{3}}\int_{-\infty}^{\infty}\frac{1}{\left|\left(\frac{2}{3}+a + iy\right)^2\right|^2}dy = c t^{-\frac{4}{3}} \end{equation}
That is
$$t^{\gamma_k}\nu^2(t)u^2(t)\leq c t^{\beta_k}, \ws\ws\ws \mbox{where} \ws\ws \beta_k < -1, $$
and therefore $S_1$ is bounded. Similarly, if instead of $x = -\frac{1}{3}$ we use $x=-a$ in (32), we obtain by (25) that $\nu^2(t)u^2(t)\leq c M^2 t^{2a-2}$, and so by (4)
\begin{equation} I = \int_0^{\frac{A}{C}}\frac{u^2(Ct)\nu^2(Ct)}{\nu^2(t)\mu^2(t)}dt\leq c \int_0^{\frac{A}{C}}\frac{t^{2(a-1)}}{\nu^2(t)} <\infty \end{equation}

This proves Theorem 1.

\medskip

We now turn to the proof of Theorem \ref{thm2}.
We will need a technical lemma. Following carefully the proof of Lemma 7 -- Lemma 11 in \cite{fu}, actually W. Fuchs proved the following:

\begin{lemma}\cite{fu}
If there is a nonnegative, monotone increasing function $h$ on $(0,\infty)$,  which fulfils (12), and
\begin{equation}\int_1^{\infty}\frac{h(r)}{r^2} = \infty,\end{equation}
and if there is a function $g$ regular on $\Re z \geq 0$ such that there are $C,c >0$, $\alpha >0$
\begin{equation}|g(z)| \leq C\left(\frac{cx}{h(r)}\right)^{\frac{x}{\alpha}}, \end{equation}
then
\begin{equation} g\equiv 0 \ws\ws\ws \mbox{on} \ws\ws \Re z \geq 0\end{equation}\end{lemma}

\remark

\medskip

In Lemma 2 $C$ and $c$ means that instead of a regular function $g$ another regular function: $bA^zg(z)$ can be considered  ($A,b$ are positive constants). It means that
$\Psi(r)$ can be  replaced by a function $\Psi_1(r)$ such that $\frac{\Psi}{\Psi_1}$ lies between finite positive bounds, and $\Psi_1(r)$ has a continuous derivative.  Therefore in the followings we will assume that $\Psi(r)$, that is $m(r)$, is continuously differentiable, if it is necessary. Furthermore since $m(r)$ is increasing, we will assume that the derivative of $m$ is nonnegative. If it is necessary, we can assume the same on $h$.

\medskip

\proof  of Theorem 2.

From (13), and the previous lemma it follows that if  a function $g(z)$ is regular on $\Re z \geq 0$, and it satisfies the inequality
\begin{equation} |g(z)| \leq \left(c\frac{\varphi(x)}{\Psi(r)}\right)^x, \end{equation}
then $g \equiv 0$. Namely, if $r \geq x >0$, then (13) and (37) together gives (35), and by the definition of $\varphi$ and $\Psi$, $\lim_{x \to 0+}\left(c\frac{\varphi(x)}{\Psi(r)}\right)^x=\|w\|_{2,(0,\infty)},$ so we can choose a constant $C$, such that $\left(\frac{\varphi(x)}{\Psi(r)}\right)^x \leq C \left(\frac{cx}{h(r)}\right)^{\frac{x}{\alpha}}$ on $\Re z \geq 0$.

Now let us  assume, by contradiction, that $S$ is not dense in $L^2_w$. In this case there is a function $f \not\equiv  0$ in $L^2_w$, such that the function
\begin{equation} G(z) = \int_0^{\infty}t^zf(t)w^2(t)dt\end{equation}
defined on $\Re z \geq 0$, satisfies the equalities
\begin{equation} G(a_k) =0 \ws\ws\ws k=1,2, \dots \end{equation}
and we can estimate its modulus by
\begin{equation}|G(z)| \leq \|f\|_{2,w}\left(\varphi(x)\right)^x \end{equation}
Let us define now on $\Re z \geq 0$
\begin{equation} g(z) = \frac{G(z)}{H(z) C_1^{z+1}}, \end{equation}
where $H$ is as in (19). By Lemma 4 \cite{fu}
\begin{equation} |H(z)| \geq \left(C_2\Psi(r)\right)^x \end{equation}
on $\mathbb{C} \setminus \cup_{k=1}^{\infty}B\left(a_k, \frac{d}{3}\right)$, where $B\left(a_k, \frac{d}{3}\right)$ are balls around $a_k$ with radius depending on $d$ (see (11)), and on the imaginary axis without exception. This implies that
$$|g(z)| \leq \left(c\frac{\varphi(x)}{\Psi(r)}\right)^x$$
on $\Re z \geq 0 \setminus \cup_{k=1}^{\infty}B\left(a_k, \frac{d}{3}\right)$ (and on the imaginary axis). But $g$ is regular on $\Re z \geq 0$, so this inequality holds on the whole half-plane, and thus by Lemma 2 $g\equiv 0$, and hence $G\equiv 0$, a contradiction.

\medskip

\proof  of Theorem \ref{thm4}.

Let us introduce the following notation on $0\leq x \leq r$, where $r\geq 0$ is fixed:
\begin{equation}v_r(x) = 2xm(r)-\frac{1}{2}\log K(x)\end{equation}

Because $\frac{K^{'}}{K}$ is increasing (see the first remark after the proof of Theorem 1), $v_r(x)$ is concave on $[0,r]$. That is, we need to distinguish three cases: (a) $v_r(x)$ is strictly decreasing on $(0,r]$, (b) $v_r(x)$ has a maximum on $(0,r]$, (c) $v_r(x)$ is strictly increasing on $[0,r]$.

\medskip

In case (a) the first derivative of $v_r(x)$ is negative on $[0,r]$, that is
\begin{equation}2m(r) < \frac{K^{'}}{2K}(x) \ws\ws\ws  0\leq x \leq r. \end{equation}
Furthermore $\frac{K^{'}}{K}$ is increasing, and it means that
\begin{equation}2m(r) \leq \frac{K^{'}}{2K}(0) \end{equation}
Since the right-hand side is constant, and the left-hand side is increasing, there is an $r_0$, such that for all $r>r_0$ (45) must be wrong.

\medskip

In case (b) there is an $0<x_0=x_0(r)\leq r$, where $v'_r(x_0)=0$. That is, for all $0\leq x \leq r$
\begin{equation}v_r(x) \leq v_r(x_0) = \frac{x_0}{2}\left(\frac{K^{'}}{K}(x_0)- \frac{\log K(x_0)}{x_0}\right)\leq \frac{r}{2}\left(\frac{K^{'}}{K}(r)- \frac{\log K(r)}{r}\right) \end{equation}
because $\frac{x}{2}\left(\frac{K^{'}}{K}(x)- \frac{\log K(x)}{x}\right)$ is increasing, since it's derivative is $\frac{1}{2}x\left(\frac{K^{'}}{K}\right)^{'}(x)$, which is nonnegative.

\medskip

In case (c) $2m(r) > \frac{K^{'}}{2K}(x)$  if $ 0\leq x \leq r$. That is, $2m(r) > \frac{K^{'}}{2K}(r)$. In this case $v_r(x) \leq v_r(r)$, $\frac{r}{2}\left(\frac{K^{'}}{K}(r)- \frac{\log K(r)}{r}\right) \leq v_r(r)$, and $v_r(r)$ itself is increasing, because using the remark after Lemma 1
$$v_r(r)^{'} = \left(2m(r)- \frac{1}{2}\frac{K^{'}}{K}(r)\right) + 2rm^{'}(r)\geq 0$$
That is, we can find a constant $C$, such that $v_r(x) \leq f(r)$ even in case (a), and $f$ is increasing.


\begin{thebibliography}{99}

\bibitem {a} J. M. Almira, M\"untz Type Theorems, {\it Surveys in Approximation Theory} {\bf 3} (2007) 152-194.
\bibitem {a1} J. M. Almira, On M\"untz Theorem for Countable Compact Sets, {\it Bull. Belg. Math. Soc. Simon Stevin} {\bf 13} 69-73
\bibitem {ba} G. V. Badalyan, On a Theorem of Fuchs, {\it Mat. Zametki} {\bf 5} (6) (1969) 723-731.
\bibitem {be} S. N. Bernstein, Sur lesrecherches r\'ecentes a la meilleure approximation des fonctions continues par les polynomes, in {\it Proc. of 5th Inter. Math. Congress} {\bf 1} (1912) 256-266.
\bibitem {bo} R. P. Boas, Density Theorems for Power Series and Complete Sets, {\it Trans. Amer. Math. Soc.} {\bf 61} (1) (1947) 54-68.
\bibitem {bopo} R. P. Boas, H. Pollard, Properties Equvivalent to the Completeness of $\{e^{-t}t^{\lambda_n}\}$, {\it Bull. Amer. Math. Soc.} {\bf 52} (4) (1946) 348-351.
\bibitem {be1} P. Borwein, T. Erd\'elyi, The Full M\"untz Theorem in $C[0,1]$ and $L_1[0,1]$, {\it J. London Math. Soc.} {\bf 54} (1996) 102-110.
\bibitem {be2} P. Borwein, T. Erd\'elyi, Generalizations of M\"untz Theorem via a Remez-type Inequality for M\"untz Spaces, {\it J. Amer. Math. Soc.} {\bf 10} (1997) 327-349.
\bibitem {be3} P. Borwein, T. Erd\'elyi, M\"untz's Theorem on Compact Subsets of Positive Measure, in {\it Approximation Theory, eds.Govil et. al., Marcel Dekker} (1988) 115-131.
\bibitem {be4} P. Borwein, T. Erd\'elyi, Polynomials and Polynomial Inequalities, {\it Springer Verlag, New York} (1996)
\bibitem {de} G. T. Deng, Incompleteness and Closure of a Linear Span of Exponential System in a Weighted Banach Space, {\it J. Approx. Theory} {\bf 125} (2003) 1-9.
\bibitem {e} T. Erd\'elyi, The "Full M\"untz Theorem" Revisited, {\it Constr. Approx.} {\bf 21} (2005) 319-335.
\bibitem {ej} T. Erd\'elyi, W. B. Johnson, The Full M\"untz Theorem in $L_p(0,1)$ for $0<p<\infty$, {\it J. Anal. Math.} {\bf 84} (2001) 145-152.
\bibitem {fu} W. H. J. Fuchs, On the Closure of $\left\{e^{-t}t^{\alpha_{\nu}}\right\}$, {\it Proc. Cambr. Phil. Soc}, {\bf 42} (1946) 91-105.
\bibitem {he} P. Henrici, Applied and Computational Complex Analysis, {\it Wiley} (1991)
\bibitem {kh} B. N. Khabibullin, Completeness of Exponential Systems and Uniqueness Sets, Russian. {\it Russia, Ufa, Bashkir State University Press} (2006)
\bibitem {ks} A. Kro\'o, J. Szabados, On Weighted Approximation by Lacunary polynomials and Rational Functions on the Half-axis, {\it East J. on Approx.} {\bf 2} (3) (1996) 289-300.
\bibitem {le} A. F. Leontev, On the Problem of the Completeness on a System of Powers on the Semi-axis, Russian. {\it Izv. Akad. Nauk SSSR, Ser. Matem.} {\bf 26} (5) (1962) 781-792.
\bibitem {lelu} E. Levin, D. S. Lubinsky, Orthogonal Polynomials for Exponential Weights, {\it Springer-Verlag New York} (2001)
\bibitem {ma} P. Malliavin, Sur qelques procedes d'extrapolation, {\it Acta Math.} {\bf 93} (1955) 179-255.
\bibitem {mu} Ch. H. M\"untz, \"Uber den Approximationssatz von Weierstrass, in {\it H. A. Schwarz's Festschrift, Berlin} (1914) 303-312.
\bibitem {o} V. Operstein, Full M\"untz Theorem in $L_p(0,1)$, {\it J. Approx. Theory} {\bf 85} (1996) 233-235.
\bibitem {p} A. Pinkus, Density in Approximation Theory, {\it Surveys in Approximation Theory} {\bf 1} (2005) 1-45.
\bibitem {sze} G. Szeg\H o, Orthogonal Polynomials, {\it Amer. Math. Soc. Pupl.} {\bf XXIII} (1959)
\bibitem {visha} B. V. Vinnitskii, A. V. Shapovalovskii, A Remark on the Completeness of Systems of Exponentials with Weight in $L_2(\mathbb{R})$, {\it Ukrainian Math. J.} {\bf 52} (7) (2000) 1002-1009.
\bibitem {wiwa} E. T. Whittaker, G. N. Whatson, Modern Analysis, {\it Cambridge Univ. Press} (1952)
\bibitem {za} R. A. Zalik, Weighted Polynomial Approximation on Unbounded Intervals, {\it J. Approx. Theory} {\bf 28} (1980) 113-119.
\bibitem {zi} E. Zikkos, Completeness of an Exponential System In Weighted Banach Spaces and Closure of its Linear Span, {\it J. Approx. Theory} {\bf 146} (2007) 115-148.
\end{thebibliography}
\end{document}